\documentclass{article}
\usepackage{spconf}
\usepackage{epsfig}
\usepackage{amsmath}
\usepackage{amsthm}
\usepackage{amsfonts}
\usepackage{amssymb}
\usepackage{mdwlist}
\usepackage{graphicx,enumerate}

\newtheorem{theorem}{Theorem}
\newtheorem{lemma}{Lemma}

\newtheorem{proposition}{Proposition}

\newtheorem{definition}{Definition}

\newcommand{\va}{{{\alpha}}}

\newcommand{\I}{\mathbf{I}}
\newcommand{\vx}{{x}}

\newcommand{\R}{\mathbb{R}}

\newcommand{\T}{\mathbf{T}}

\newcommand{\be}{\begin{eqnarray}}
\newcommand{\ee}{\end{eqnarray}}

\newcommand{\M}{\mathcal{M}}
\renewcommand{\P}{\mathsf{P}}

\newcommand{\C}{\mathcal{C}}

\newcommand{\norm}[1]{\left\|#1\right\|}
\newcommand{\abs}[1]{\left|#1\right|}

\newcommand{\parenth}[1]{\left({#1}\right)}

\newcommand{\eps}{\varepsilon}
\newcommand{\Hm}{\mathcal{H}}
\newcommand{\Tr}{^\mathrm{T}}

\renewcommand{\ge}{\geqslant}
\renewcommand{\le}{\leqslant}
\renewcommand{\geq}{\geqslant}
\renewcommand{\leq}{\leqslant}
\newcommand{\reg}{f_2}
\newcommand{\Prj}{\mathcal{P}}

\DeclareMathOperator{\rprox}{rprox}
\DeclareMathOperator{\prox}{prox}
\DeclareMathOperator*{\argmin}{arg\ min}

\title{Image Deconvolution under Poisson Noise using Sparse Representations
  and Proximal Thresholding Iteration}
\name{F.-X. Dup\'e$^\text{a}$, M.J. Fadili$^\text{a}$ and J.-L. Starck$^\text{b}$}
\address{\begin{tabular}{cc}
        $^\text{a}$ GREYC UMR CNRS 6072 & $^\text{b}$ DAPNIA/SEDI-SAP CEA-Saclay\\
        14050 Caen France & 91191 Gif-sur-Yvette France
        \end{tabular}
}

\begin{document}
%\ninept
%
\small
\maketitle
\begin{abstract}
We propose an image deconvolution algorithm when the data is contaminated by Poisson noise. 
The image to restore is assumed to be sparsely represented in a dictionary of waveforms such as the wavelet or curvelet transform. Our key innovations are: First, we handle the Poisson noise properly by using the Anscombe variance stabilizing transform leading to a {\it non-linear} degradation equation with additive Gaussian noise. Second, the deconvolution problem is formulated as the minimization of a convex functional with a data-fidelity term reflecting the noise properties, and a non-smooth sparsity-promoting penalties over the image representation coefficients (e.g. $\ell_1$-norm). Third, a fast iterative backward-forward splitting algorithm is proposed to solve the minimization problem. We derive existence and uniqueness conditions of the solution, and establish convergence of the iterative algorithm. Experimental results are carried out to show the striking benefits gained from taking into account the Poisson statistics of the noise. These results also suggest that using sparse-domain regularization may be tractable in many deconvolution applications, e.g. astronomy or microscopy.
\end{abstract}
\section{Introduction}
\label{sec:intro}
Deconvolution is a longstanding problem in many areas of signal and image processing (e.g. biomedical imaging \cite{Sarder2006}, astronomy \cite{Starck2006}, remote-sensing, to quote a few). For example, research in astronomical image deconvolution has recently seen considerable work, partly triggered by the Hubble space telescope (HST) optical aberration problem at the beginning of its mission. In biomedical imaging, researchers are also increasingly relying on deconvolution to improve the quality of images acquired by confocal microscopes. Deconvolution may then prove crucial for exploiting images and extracting scientific content.

There is an extensive literature on deconvolution problems. One might refer to well-known dedicated monographs on the subject. In presence of Poisson noise, several deconvolution methods have been proposed such as Tikhonov-Miller inverse filter and Richardson-Lucy (RL) algorithms; see \cite{Sarder2006,Starck2006} for an excellent review. The RL has been used extensively in applications because it is adapted to Poisson noise. The RL algorithm, however, amplifies noise after a few iterations, which can be avoided by introducing regularization. In \cite{Dey2004}, the authors presented a Total Variation (TV)-regularized RL algorithm, and Starck et al. advocated a wavelet-regularized RL algorithm \cite{Starck2006}.

In the context of deconvolution with gaussian white noise, sparsity-promoting regularization over orthogonal wavelet coefficients has been recently proposed \cite{Daubechies2004,Combettes2005}. Generalization to frames was proposed in \cite{Teschke2007,Combettes2007b}. In \cite{Fadili2006a}, the authors presented an image deconvolution algorithm by iterative thresholding in an overcomplete dictionary of transforms. However, all sparsity-based approaches published so far have mainly focused on Gaussian noise.

In this paper, we propose an image deconvolution algorithm for data blurred and contaminated by Poisson noise. The Poisson noise is handled properly by using the Anscombe VST, leading to a {\it non-linear} degradation equation with additive Gaussian noise, see \eqref{eq:3}. To regularize the solution, we impose a sparsity prior on the representation coefficients of the image to restore, e.g. wavelet or curvelet coefficients. Then, the deconvolution problem is formulated as the minimization of a convex functional with a data-fidelity term reflecting the noise properties, and a non-smooth sparsity-promoting penalties over the image representation coefficients. Inspired by the work in \cite{Combettes2005}, a fast proximal iterative algorithm is proposed to solve the minimization problem. We also provide an analysis of the optimization problem and establish convergence of the iterative algorithm. Experimental results are carried out to compare our approach and show the striking benefits gained from taking into account the Poisson nature of the noise.

\subsection*{Notation}
\label{sec:notation}
Let $\mathcal{H}$ a real Hilbert space, here a finite
dimensional vector subspace of $\R^n$. We denote by
$\norm{.}_2$ the norm associated with the inner product in
$\mathcal{H}$, and $\I$ is the identity operator on
$\mathcal{H}$. $\vx$ and $\va$ are respectively reordered vectors of
image samples and transform coefficients. A function $f$ is proper if it 
is not identically $+\infty$ everywhere. A function $f$ is coercive, 
if $\lim_{\norm{\vx}_2 \rightarrow +\infty}f\parenth{\vx}=+\infty$.
$\Gamma_0(\mathcal{H})$ is the class of all proper lower semi-continuous convex
functions from $\mathcal{H}$ to $]-\infty,+\infty]$. The subdifferential of $f$
is denoted $\partial f$.

%The subdifferential of $f$
%is the set-valued operator $\partial f : \mathcal{H} \to 2^{\mathcal{H}}$
%whose value at $x\in\mathcal{H}$ is:
%\begin{equation*}
%  \partial f(x) = \left\{
%    u \in \mathcal{H}\ |\ (\forall y \in \mathcal{H})
%    \psp{y-x}{u} + f(x) f(y)
%  \right\}
%\end{equation*}

\section{Problem statement}
Consider the image formation model where an input image $\vx$ is blurred by a point spread function (PSF) $h$ and contaminated by Poisson noise. The observed image is then a discrete collection of counts $y=(y_i)_{1 \leq i \leq n}$ where $n$ is the number of samples. Each count $y_i$ is a realization of an independent Poisson random variable with a mean $(h \circledast x)_i$, where $\circledast$ is the circular convolution operator. Formally, this writes $y_i \sim \mathcal{P}\parenth{(h \circledast x)_i}$.

A naive solution to this deconvolution problem would be to apply traditional approaches designed for Gaussian noise. But this would be awkward as (i) the noise tends to Gaussian only for large mean $(h \circledast x)_i$ (central limit theorem), and (ii) the noise variance depends on the mean anyway. A more adapted way would be to adopt a bayesian framework with an adapted anti-log-likelihood score reflecting the Poisson statistics of the noise. Unfortunately, doing so, we would end-up with a functional which does not satisfy some key properties (the Lipschitzian property in \cite{Combettes2005}), hence preventing us from using the backward-forward splitting proximal algorithm to solve the optimization problem. To circumvent this difficulty, we propose to handle the noise statistical properties by using the Anscombe VST defined as
\begin{equation}
  \label{eq:2}
  z_i = 2\sqrt{y_i + \tfrac{3}{8}}, ~ 1 \leq i \leq n .
\end{equation}
Some previous authors \cite{ChauxSPIE} have already suggested to use the Anscombe VST, and then deconvolve with wavelet-domain regularization as if the stabilized observation $z_i$ were linearly degraded by $h$ and contaminated by additive Gaussian noise. But this is not valid as standard asymptotic results of the Anscombe VST state that
\begin{equation}
  \label{eq:3}
  z_i = 2\sqrt{(h\circledast x)_i+\tfrac{3}{8}} + \eps,\quad \eps \sim \mathcal{N}(0,1)
\end{equation}
where $\eps$ is an additive white Gaussian noise of unit variance. In words, $z$ is {\em non-linearly} related to $x$. In Section~\ref{sec:sparse-iter-thresh}, we provide an elegant optimization problem and a fixed point algorithm taking into account such a non-linearity.

\section{Sparse image representation}
\label{sec:sparse-image-repr}
Let $x \in \Hm$ be an $\sqrt{n}\times\sqrt{n}$ image. $x$ can be written as the superposition of elementary atoms $\varphi_\gamma$ parameterized by $\gamma \in \mathcal{I}$ such that:
\begin{equation}
  \label{eq:4}
  x = \sum_{\gamma \in \mathcal{I}} \alpha_\gamma \varphi_\gamma = \Phi \va,\quad \abs{\mathcal{I}} = L
\end{equation}
We denote by $\Phi$ the dictionary i.e. the $n\times L$ matrix whose
columns are the generating waveforms $\parenth{\varphi_\gamma}_{\gamma \in \mathcal{I}}$ all normalized to a unit $\ell_2$-norm.
The forward transform is then defined by a non-necessarily square matrix
$\T = \Phi\Tr \in \mathbb{R}^{L\times n}$ with $L\ge n$. When $L > n$ the dictionary is said
to be redundant or overcomplete. In the case of the simple orthogonal basis, the inverse transform is trivially $\Phi
= \T\Tr$. Whereas assuming that $\T$ is a tight frame implies that the
frame operator satisfies $\T\Tr\T = A\I$, where $A > 0$ is the tight
frame constant. For tight frames, the pseudo-inverse reconstruction operator reduces to $A^{-1}\T$.

Our prior is that we are seeking for a good representation of
$x$ with only few significant coefficients. This makes sense since most practical signals or images are compressible in some transform domain (e.g. wavelets, curvelets, DCT, etc). These transforms generally correspond to an orthogonal basis or a tight frame. In the rest of the paper, $\Phi$ will be an orthobasis or a tigth frame of $\Hm$.

\section{Sparse Iterative Deconvolution}
\label{sec:sparse-iter-deconv}

We first define the notion of a proximity operator, which was introduced in \cite{Moreau1962} as a generalization of the notion of a convex projection operator.
\begin{definition}[Moreau\cite{Moreau1962}]
  \label{def:1}
  Let $\varphi \in \Gamma_{0}(\Hm)$. Then, for every $x\in\Hm$, the function
  $y \mapsto \varphi(y) + \norm{x-y}^{2}/2$ achieves its infimum at a unique point denoted by
  $\prox_{\varphi}x$. The operator $\prox_{\varphi} : \Hm \to \Hm$ thus defined is
  the \textit{proximity operator} of $\varphi$. 
  We define the reflection operator $\rprox_{\varphi} = 2\prox_{\varphi} - \I$.
%Moreover, $\forall x,p \in \Hm$
%\begin{align*}
%  p = \prox_{\varphi} x \iff & x-p \in \partial\varphi(p)\\
%  \iff &  <\!y-p|x-p\!> + \varphi(p) \le \varphi(y)\\
%  & (\forall y\in\Hm)
%\end{align*}
\end{definition}

\vspace*{-0.5cm}
\subsection{Optimization problem}
\label{sec:sparse-iter-thresh}

The class of minimization problems we are interested in can be stated in the general form \cite{Combettes2005}:
\begin{equation}
  \label{eq:6}
  \argmin_{x \in \Hm} f_1(x) + f_2(x) .
\end{equation}
where $f_1 \in \Gamma_0(\mathcal{H})$, $f_2\in\Gamma_0(\mathcal{H})$
and $f_1$ is differentiable with $\kappa$-Lipschitz
gradient. We denote by $\M$ the set of solutions.

\noindent
From \eqref{eq:3} and \eqref{eq:4}, we immediately deduce the data fidelity term
\begin{gather}
\label{eq:7}
  F \circ H \circ \Phi~(\va), ~ \text{with} \\
  F : \eta \mapsto \sum_{i=1}^n f(\eta_i),\quad f(\eta_i) = \frac{1}{2}
  \left( z_i - 2\sqrt{\eta_i+\tfrac{3}{8}} \right)^2 , \nonumber 
\end{gather}
where $H$ denotes the (block-Toeplitz) convolution operator. From a statistical perspective, \eqref{eq:7} corresponds to the anti-log-likelihood score. 

Adopting a bayesian framework and using a standard maximum a posteriori (MAP) rule, our goal is to minimize the following functional with respect to the representation coefficients $\alpha$ 
\begin{gather}
\label{eq:9}
(\P_{\lambda,\psi}): \argmin_{\alpha} J(\alpha) \\
J : \alpha \mapsto \underbrace{F \circ H \circ \Phi\ (\alpha)}_{f_1(\va)} + 
\underbrace{\imath_{\C} \circ \Phi\ (\alpha) + \lambda \sum_{i=1}^L \psi(\alpha_i)}_{f_2(\va)}, \nonumber 
\end{gather}
where we implicitly assumed that $(\alpha_i)_{1 \leq i \leq L}$ are independent and identically distributed with a Gibbsian density $\propto e^{-\lambda\psi(\alpha_i)}$. 
Notice that $f_2$ is not smooth.
The  penalty function $\psi$ is chosen to enforce sparsity,
$\lambda > 0$ is a regularization parameter and $\imath_{\C}$ is the indicator function of a convex set $\C$. In our case, $\C$ is the positive orthant. We remind that the positivity constraint is because we are fitting Poisson intensities, which are positive by nature. We have the following,
\begin{proposition}
\label{prop:1}
  {~} \\\vspace{-0.5cm}
  \begin{itemize*}
  \item $f_1$ is convex function. It is strictly convex if $\Phi$ is an ortho-basis and $\mathrm{ker}\parenth{H}=\emptyset$ (i.e. the spectrum of the PSF does not vanish).
  \item $f_1$ is continuously differentiable with a $\kappa$-Lipschitz gradient where
    \begin{equation}
      \label{eq:22}
      \kappa \le \parenth{\tfrac{2}{3}}^{3/2} 4 A \norm{H}_2^2 \norm{z}_{\infty} < +\infty.
    \end{equation}
  \item $(\P_{\lambda,\psi})$ is a particular case of problem \eqref{eq:6}.
  \end{itemize*}
\end{proposition}

\subsubsection{Characterization of the solution}
\label{sec:existence-solution}

\begin{proposition}
\label{prop:2}
  Since $J$ is coercive and convex, the following holds:
  \begin{enumerate*}
  \item Existence: $(\P_{\lambda,\psi})$ has at least one solution, i.e. $\M \neq \emptyset$.
  \item Uniqueness: $(\P_{\lambda,\psi})$ has a unique solution if
    $\Phi$ is a basis and $\mathrm{ker}\parenth{H}=\emptyset$, or if $\psi$ is strictly convex.
  \end{enumerate*}
\end{proposition}
\begin{proof}
  The existence is obvious because $J$ is coercive. If $\Phi$ is an ortho-basis
  and $\mathrm{ker}\parenth{H} = \emptyset$ then $f_1$ is strictly convex and so is $J$ leading
  to a strict minimum. Similarly, if $\psi$ is strictly convex, so is $f_2$, hence $J$.
  
\end{proof}

\subsubsection{Proximal iteration}
\label{sec:proximal-iteration}
For notational simplicity, we denote by $\Psi$ the function $\alpha
\mapsto \sum_i \psi(\alpha_i)$. The following useful lemmas are first stated:
\begin{lemma}
\label{lemma:1}
The gradient of $\nabla f_1$ is
\begin{equation}
  \label{eq:19}
  \nabla f_1 (\va) = \Phi\Tr \circ H^* \circ \nabla F \circ H \circ \Phi\ (\va)
\end{equation}
with
\begin{equation}
  \label{eq:20}
  \nabla F (\eta) = \parenth{\frac{-z_i}{\sqrt{\eta_i + 3/8}} + 2}_{1\le i\le n}
\end{equation}
\end{lemma}
The proof is straightforward.

\begin{lemma}
\label{lemma:2}
  Let $\Phi$ an orthobasis or a tight frame with constant $c$. 
  \begin{enumerate}
  \item If $\va \in \C'$ then $\prox_{f_2}(\va) = \prox_{\lambda\Psi}(\va)$, $\C'=\{\va | \Phi\va \in \C\}$.
  \item Otherwise, let $\sum_t \nu_t(1-\nu_t)=+\infty$. Take $\gamma^{0} \in \Hm$, and
    define the sequence of iterates:
    \begin{eqnarray}
      \label{eq:proxtframe1}
      \begin{split}
        \gamma^{t+1} = \gamma^t + \nu_t\left(\rprox_{\lambda\Psi + \tfrac{1}{2}\norm{. - \va}^2}\circ\right.\\
          \left. \rprox_{\imath_{\C'}} - \I \right) (\gamma^t) ,
      \end{split}
    \end{eqnarray}
    where\\ $\prox_{\lambda\Psi + \tfrac{1}{2}\norm{. - \va}^2}
    (\gamma^t)=\parenth{\prox_{\tfrac{\lambda}{2}\psi} ((\va_i + \gamma^t_i)/2)}_{1\leq i \leq
      L}$, $\Prj_{\C'} = \prox_{\imath_{\C'}} = c^{-1} \Phi\Tr\circ\Prj_\C\circ\Phi
    +\parenth{\I - c^{-1} \Phi\Tr\circ\Phi}$ and $\Prj_\C$ is the projector onto the positive
    orthant $(\Prj_\C \eta)_i = \max(\eta_i, 0)$. Then,
    \begin{equation}
      \label{eq:proxtframe2}
      \gamma^t \rightharpoonup \gamma ~ \text{and} ~ \prox_{\reg}(\va) = \Prj_{\C'}(\gamma).
    \end{equation}
  \end{enumerate}
\end{lemma}

\noindent
We are now ready to state our main proximal iterative algorithm to solve the minimization problem $(\P_{\lambda,\psi})$:
\begin{theorem}
  \label{th:2}
%For $t \geq 0$, let $(\mu_{t} )_{t}$ be a sequence in $]0, +\infty[$
%such that $0 < \inf_{t} \mu_{t} \le \sup_{t} \mu_{t}< \parenth{\frac{3}{2}}^{3/2}/\parenth{2 A \norm{H}_2^2 \norm{z}_{\infty}}$, let $(\theta_{t})_{t}$ be a sequence in $]0,1]$ such that $\inf_{t}\theta_{t} 0$, and let $(a_{t})_{t}$ and $(b_{t})_{t}$ be sequences in $\Hm$ such that $\sum_{t}\norm{b_{t}} < +\infty$ and similarly for $a_t$. Fix $\va_{0}\in\Hm$, for every $t \ge 0$, set
%\begin{equation}
%  \label{eq:5}
%  \va_{t+1} = \va_{t} + \theta_{t} \parenth{\prox_{\mu_{t}f_{2}}\parenth{\va_{t} -
%  \mu_{t}\parenth{ \nabla f_{1}(\va_{t}}+ b_{t})} + a_{t} - \va_{t}}
%\end{equation}
%where $\nabla f_{1}$ and $\prox_{\mu_{t}f_{2}}$ are given by Lemma \ref{lemma:1} and \ref{lemma:2}.
%\noindent 
%Then $(\alpha_t)_{t \geq 0}$ converges (weakly) to a solution of $(\P_{\lambda,\psi})$.
%\end{theorem}
For $t \geq 0$, let $(\mu_{t} )_{t}$ be a sequence in $]0, +\infty[$
such that $0 < \inf_{t} \mu_{t} \le \sup_{t} \mu_{t}< \parenth{\frac{3}{2}}^{3/2}/\parenth{2 A \norm{H}_2^2 \norm{z}_{\infty}}$. 
Fix $\va_{0}\in\Hm$, for every $t \ge 0$, set
\begin{equation}
  \label{eq:5}
  \va_{t+1} = \prox_{\mu_{t}f_{2}}\parenth{\va_{t} -
  \mu_{t}\parenth{ \nabla f_{1}(\va_{t}})}
\end{equation}
where $\nabla f_{1}$ and $\prox_{\mu_{t}f_{2}}$ are given by Lemma \ref{lemma:1} and \ref{lemma:2}.
\noindent 
Then $(\alpha_t)_{t \geq 0}$ converges (weakly) to a solution of $(\P_{\lambda,\psi})$.
\end{theorem}

\begin{proof}
We give a sketch of the proof. The main theorem on the proximal iteration can be found in
\cite[Theorem 3.4]{Combettes2005}. Hence, combining this theorem with Lemma \ref{lemma:1}, Lemma \ref{lemma:2} and Proposition \ref{prop:1}, the result follows.

\end{proof}
Note that if the PSF $h$ is low-pass normalized to a unit sum, then $\norm{H}_2^2 = 1$.

\noindent 
We now turn to $\prox_{\delta\psi}$ which is given by the following result:
\begin{theorem}
\label{th:3}
Suppose that $\psi$ satisfies, (i) $\psi$ is convex even-symmetric , non-negative and non-decreasing on $[0,+\infty)$, and $\psi(0)=0$. (ii) $\psi$ is twice differentiable on $\mathbb{R}\setminus \{0\}$. (iii) $\psi$ is continuous on $\mathbb{R}$, it is not necessarily smooth at zero and admits a positive right derivative at zero $\psi^{'}_+(0) = \lim_{h\to 0^+} \frac{\psi(h)}{h} > 0$. Then, the proximity operator $\prox_{\delta\psi}(\beta) = \hat{\va}(\beta)$ has exactly one continuous solution decoupled
in each coordinate $\beta_i$ :
\begin{equation}
\label{eq:10}
\hat{\va}_i(\beta_i) =
\begin{cases}
  0 & \text{if } \abs{\beta_i} \le \delta\psi^{'}_+(0)\\
  \beta_i-\delta\psi^{'}(\hat{\va}_i) & \text{if } \abs{\beta_i} > \delta\psi^{'}_+(0)
\end{cases}
\end{equation}
\end{theorem}
A proof of this theorem can be found in \cite{Fadili2006}. Among the most popular penalty functions $\psi$ satisfying the above requirements, we have $\psi(\va_i) = \abs{\va_i}$, in which case the associated proximity operator is soft-thresholding. In this case, iteration \eqref{eq:5} is essentially an iterative thresholding algorithm with a positivity constraint. 

\section{Experimental results}
\label{sec:results}
The performance of our approach has been assessed on several 2D datasets, from which we here illustrate two examples. Our algorithm was compared to RL without regularization, RL with multi-resolution support wavelet regularization \cite[RL-MRS]{Starck2006}, the naive proximal method that would assume the noise to be Gaussian (NaiveGauss), and the approach of \cite{ChauxSPIE} (AnsGauss). For all results presented, each algorithm was run with 200 iterations, except RL which was stopped when its MSE was smallest.

\begin{figure}[ht]
\begin{tabular}{@{}c@{}c@{}c@{}}
\includegraphics[width=0.33\linewidth]{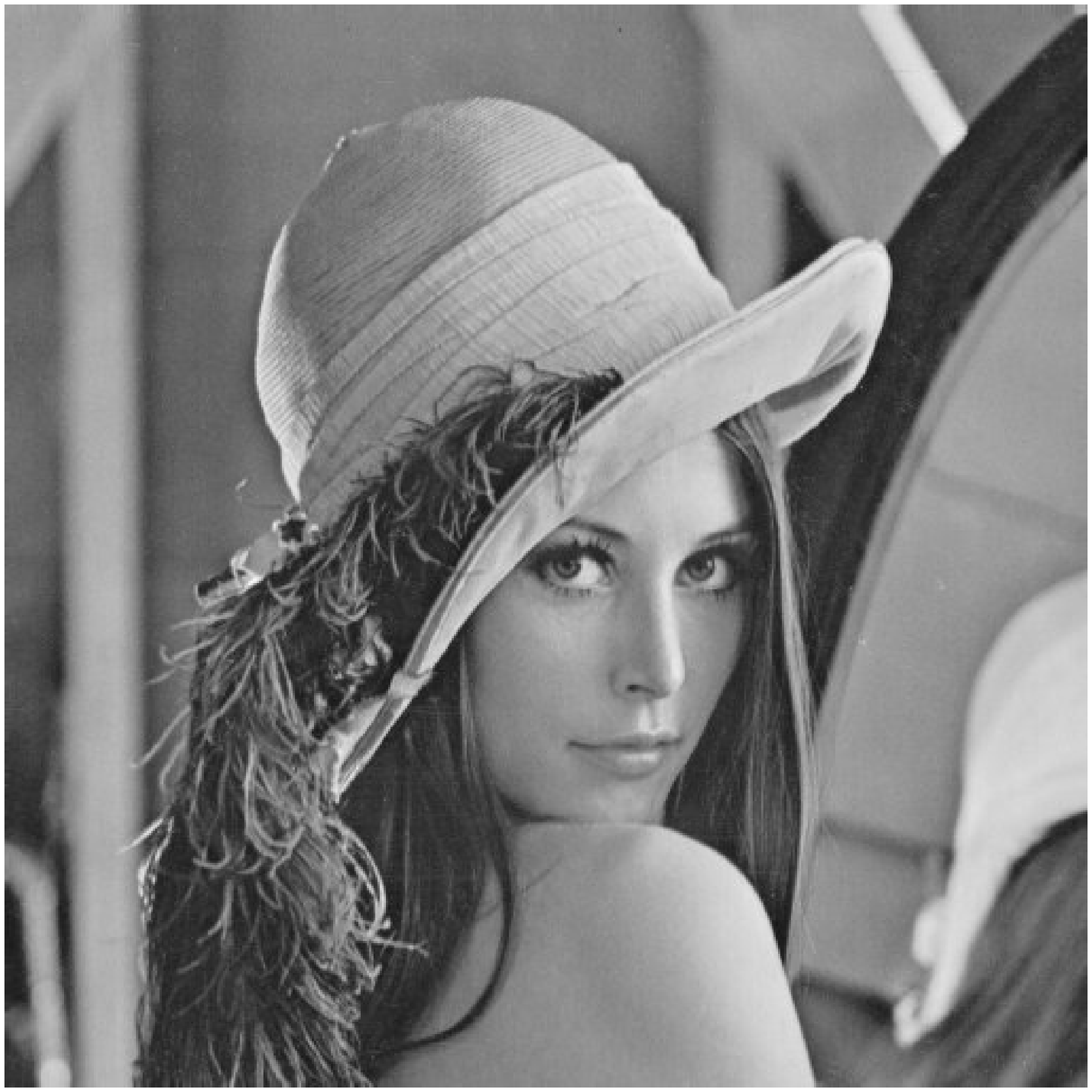} &
\includegraphics[width=0.33\linewidth]{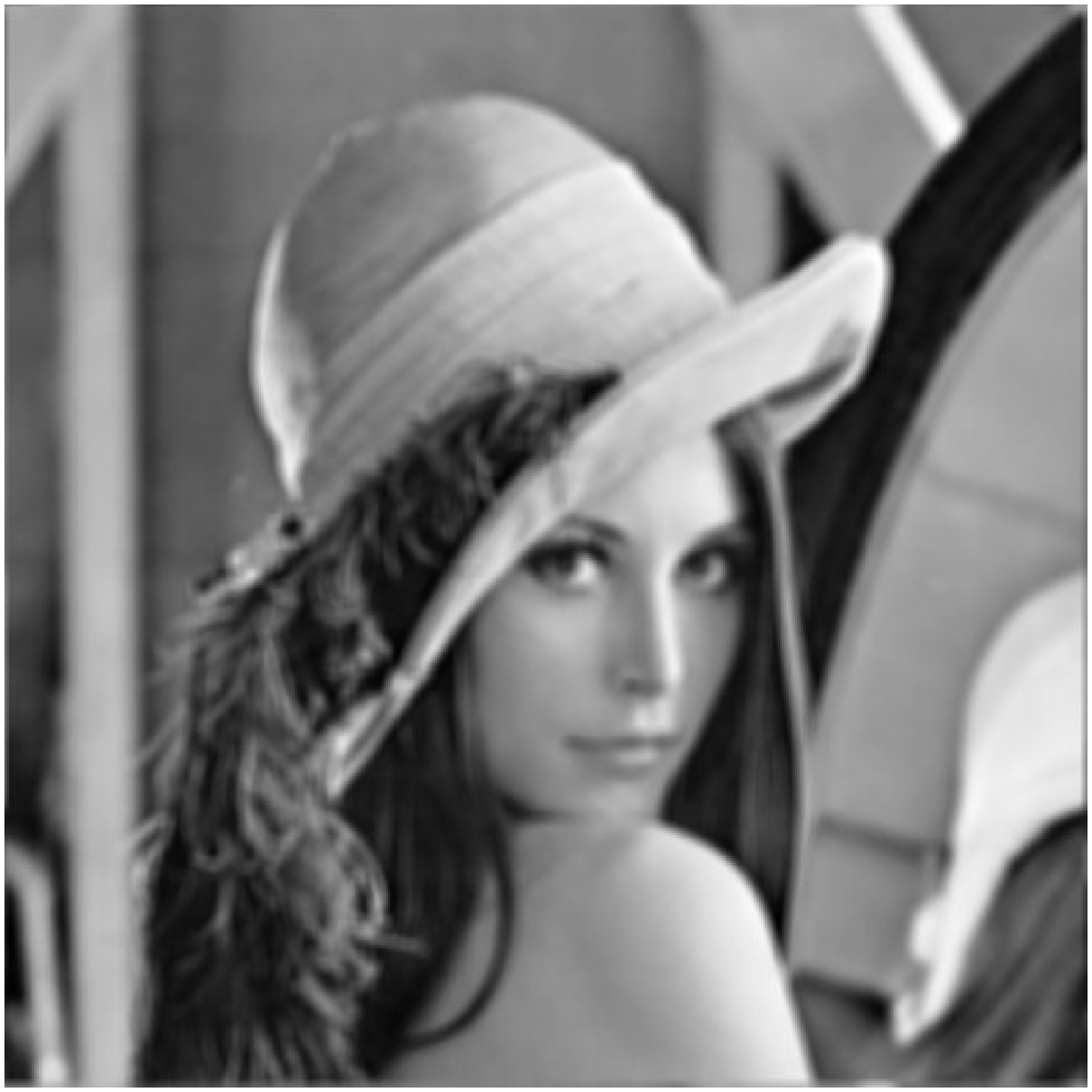} &
\includegraphics[width=0.33\linewidth]{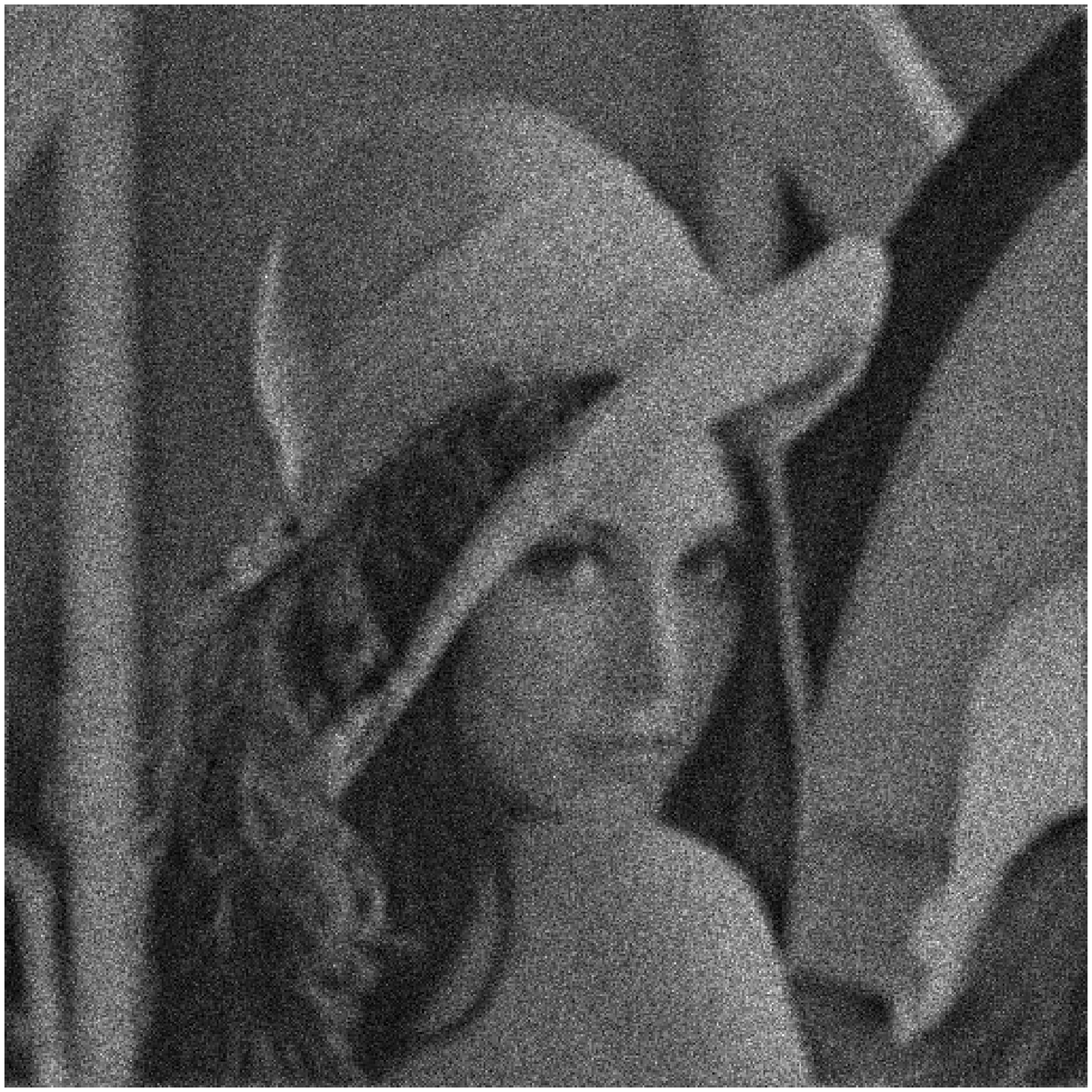} \\
(a) & (b) & (c)	\\
\includegraphics[width=0.33\linewidth]{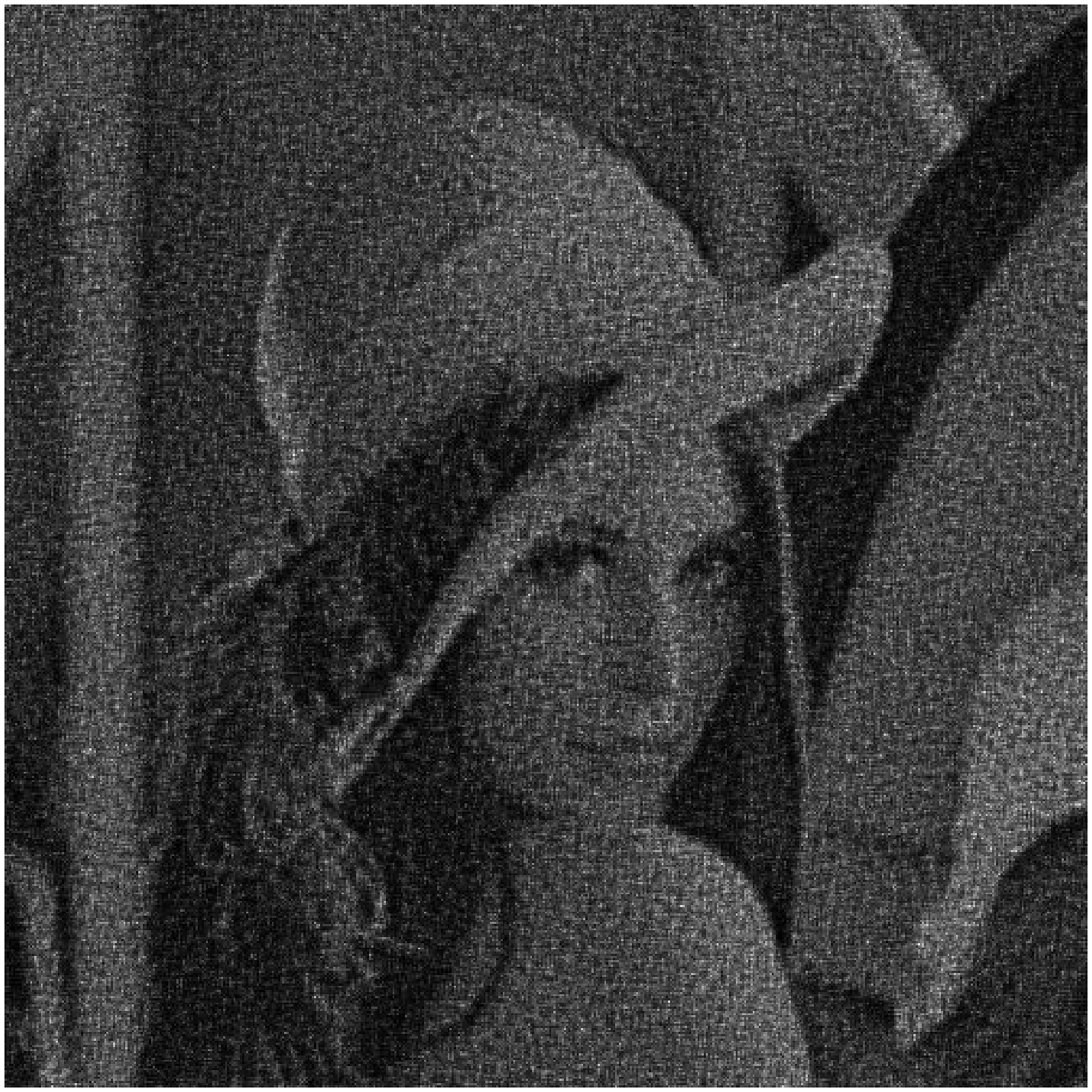} &
\includegraphics[width=0.33\linewidth]{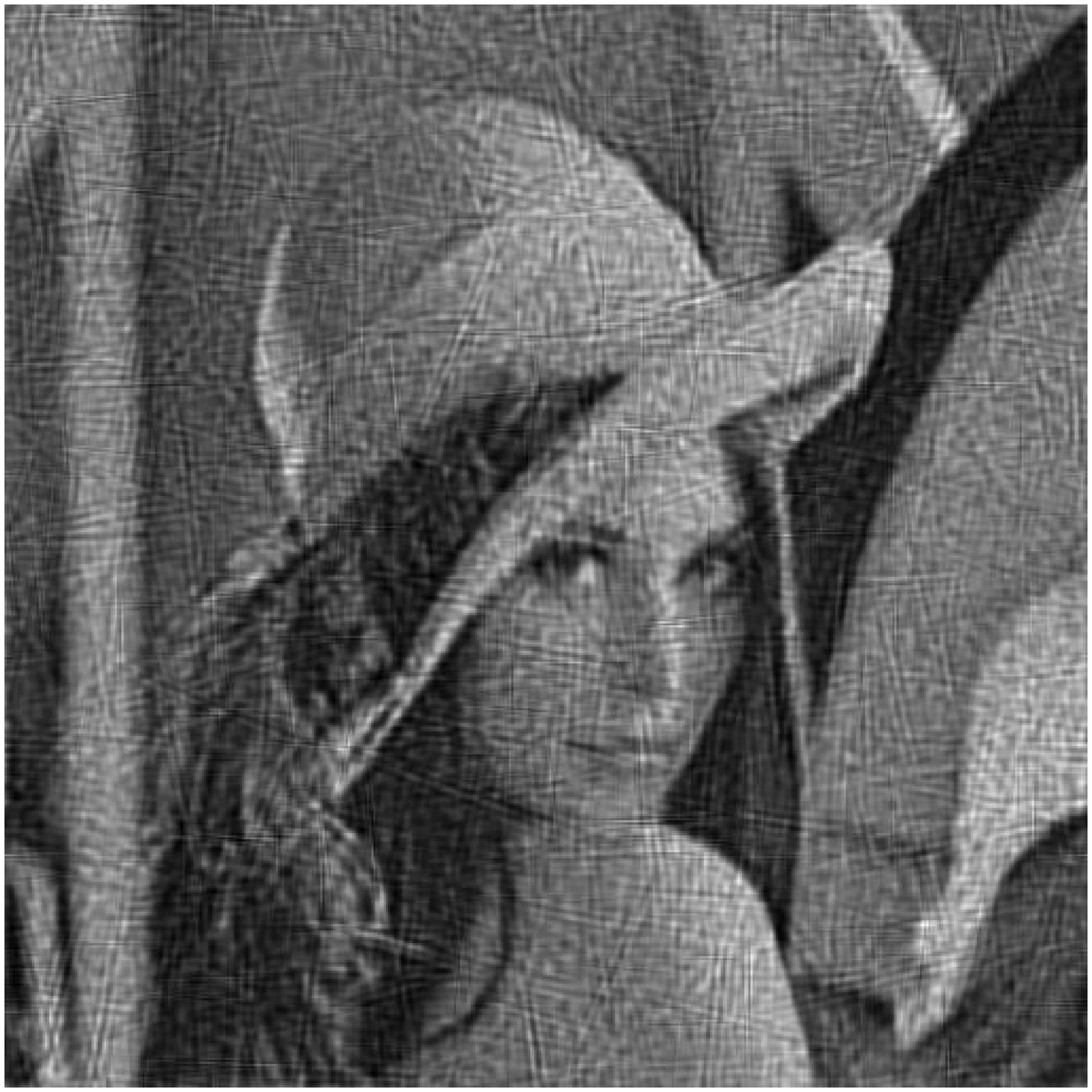} &
\includegraphics[width=0.33\linewidth]{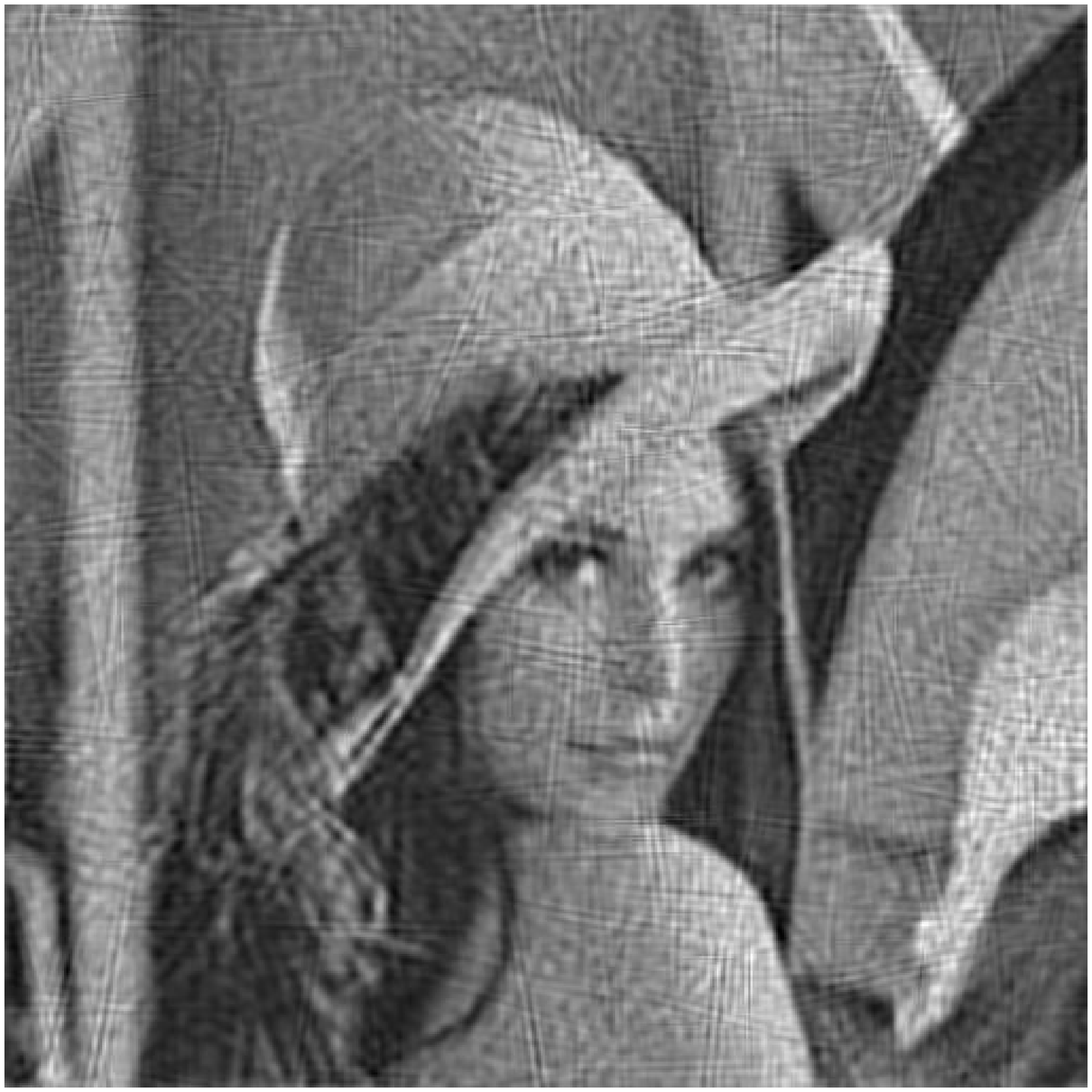} \\
(d) & (e) & (f)	\\
\includegraphics[width=0.33\linewidth]{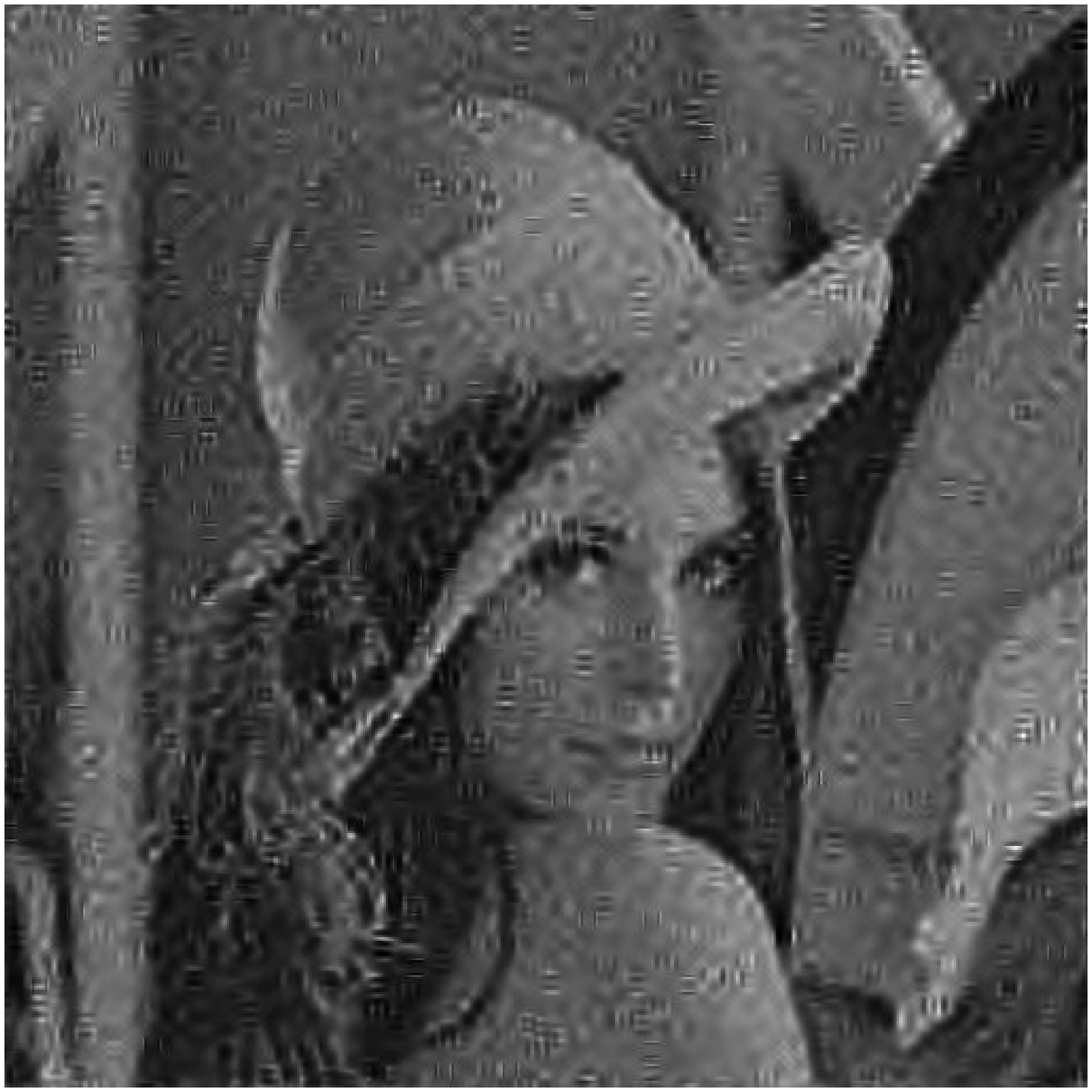} &
\includegraphics[width=0.33\linewidth]{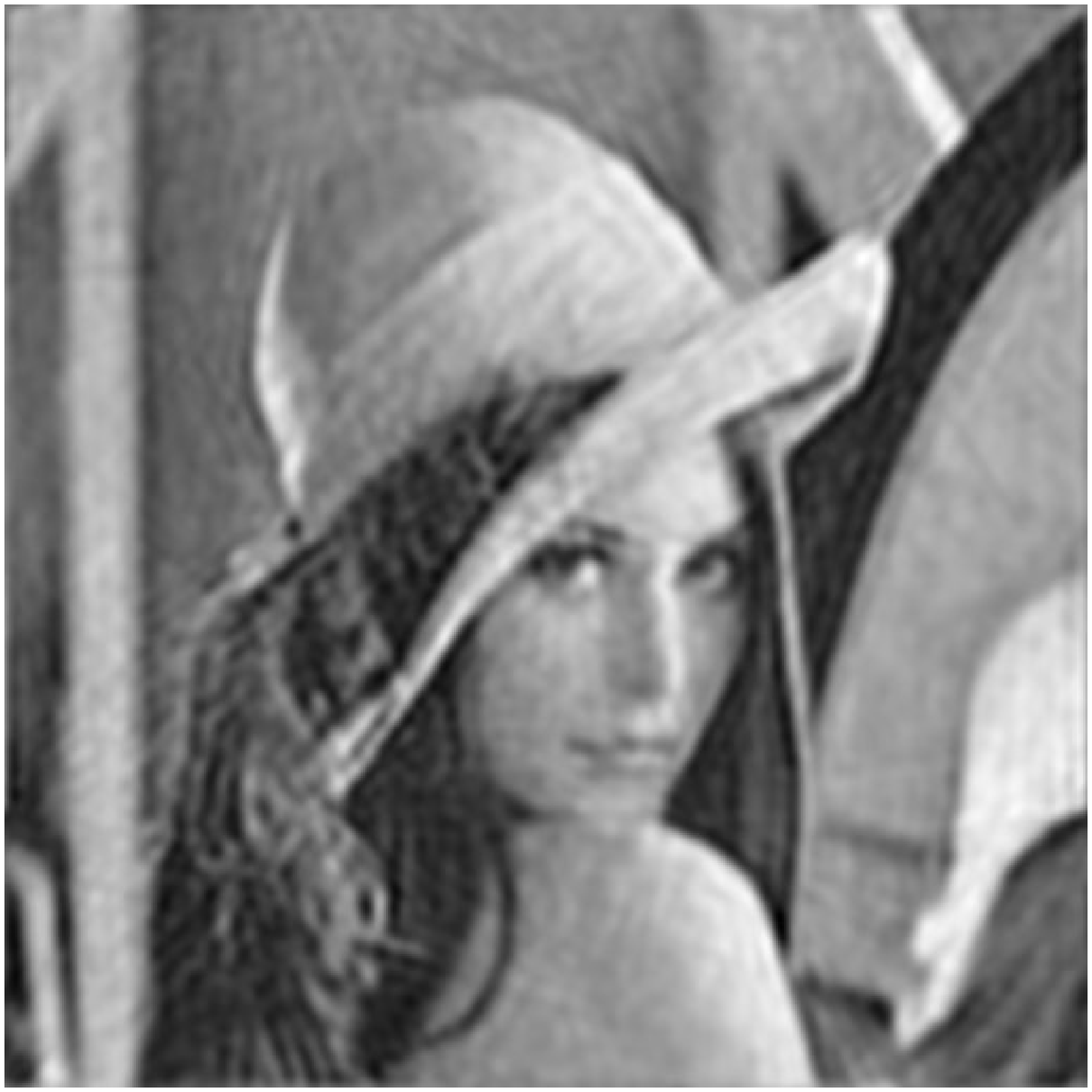} & \\
(g) & (h) &   	
\end{tabular}
\caption{Deconvolution of Lena. (a) Original, (b) Blurred, (c) Blurred+noisy, (d) RL, (e) NaiveGauss, (f) AnsGauss, (g) RL-MRS, (h) Our algorithm.}
\label{fig:lena}
\end{figure}

In Fig.\ref{fig:lena}, the original Lena image with a maximum intensity of 30 is depicted in (a), its blurred and blurred+noisy versions are in (b) and (c). With Lena, and for NaiveGauss, AnsGauss and our approach, the dictionary $\Phi$ contained the curvelet tight frame. The deconvolution results are shown in Fig.\ref{fig:lena}(d)-(h). As expected, the RL is the worst as it lacks regularization. There are also noticeable artifacts in NaiveGauss, AnsGauss and RL-MRS. Our deconvolved image appears much cleaner. This visual impression is confirmed by quantitative measures of the quality of deconvolution, where we used both the mean $\ell_1$-error (adapted to Poisson noise), and the well-known MSE criteria. The mean $\ell_1$-error for Lena is shown in Tab.\ref{tab:intens} (similar results were obtained for the MSE). In general, our approach performs very well. At low intensity levels, RL-MRS has the smallest error very comparable to our approach. For the other intensity levels, our algorithm provides the best performance. At high intensity levels, NaiveGauss is competitive. This comes as no surprise since this is an intensity regime where Poisson noise approaches the Gaussian behavior. On the other hand, the results reveal that AnsGauss performs poorly just after RL, probably because it does not handle properly the non-linearity of the degradation model \eqref{eq:3} after the VST. 

\begin{figure}[ht]
\begin{tabular}{@{}c@{}c@{}c@{}}
\includegraphics[width=0.33\linewidth]{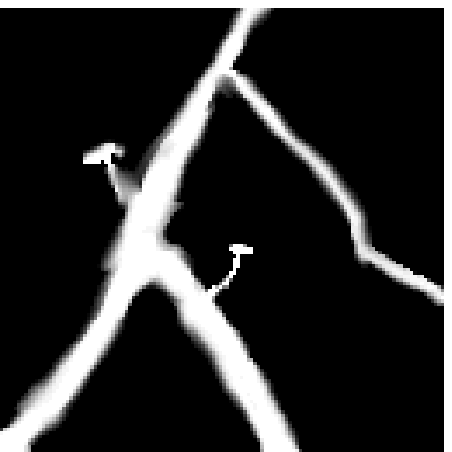} &
\includegraphics[width=0.33\linewidth]{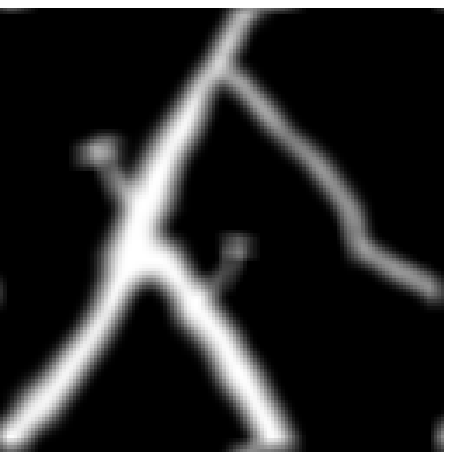} &
\includegraphics[width=0.33\linewidth]{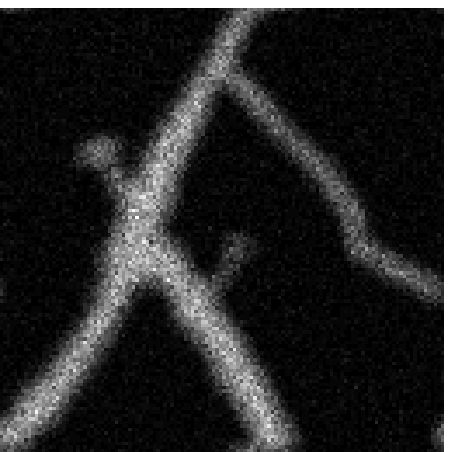} \\
(a) & (b) & (c)	\\
\includegraphics[width=0.33\linewidth]{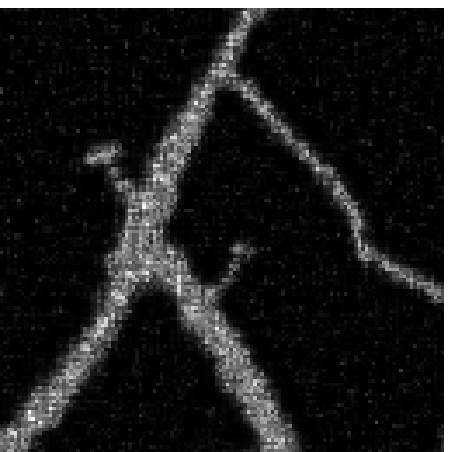} &
\includegraphics[width=0.33\linewidth]{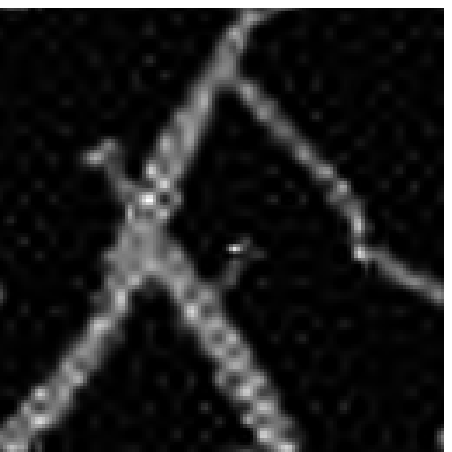} &
\includegraphics[width=0.33\linewidth]{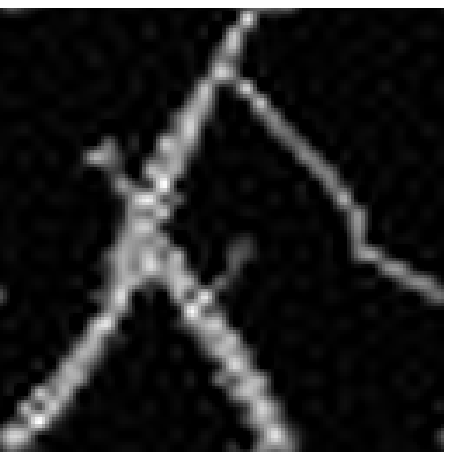} \\
(d) & (e) & (f)	\\
\includegraphics[width=0.33\linewidth]{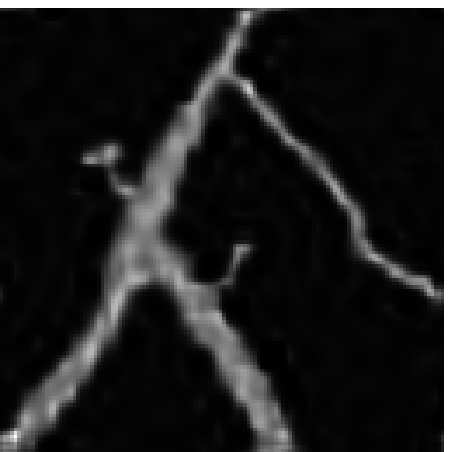} &
\includegraphics[width=0.33\linewidth]{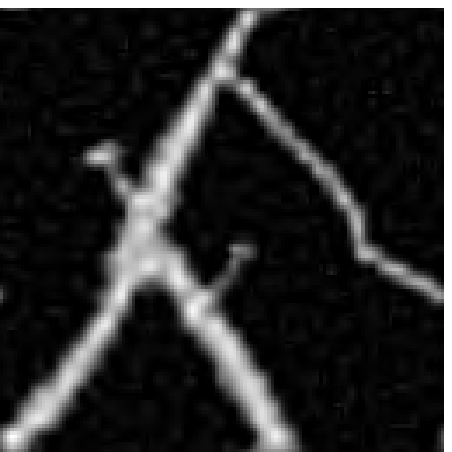} & \\
(g) & (h) &   	
\end{tabular}
\caption{Deconvolution of Neuron image. (a) Original, (b) Blurred, (c) Blurred+noisy, (d) RL (mean $\ell_1$-error 2.1841), (e) NaiveGauss (1.6368), (f) AnsGauss (2.0477), (g) RL-MRS (1.1926), (h) Our algorithm (1.335).   }
\label{fig:neuron}
\end{figure}

We further illustrate the capabilities of our approach on a confocal microscopy simulation. We have created a phantom of an image of a neuron dendrite segment with a mushroom-shaped spines, see Fig.\ref{fig:neuron}. The experimental settings were the same as for Lena except that the dictionary here contained the wavelet transform. The findings are similar to those of Lena both visually and quantitatively.

\begin{table}[ht]
  \centering
  % M1 = with good stabilization
  % M2 = without stabilizing
  % M3 = bad stabilized
  % M4 = RL-MRS
  % M5 = RL
  \begin{tabular}{|c|c|c|c|c|}
    \hline  & \multicolumn{4}{|c|}{Intensity regime}\\
    \hline Method		        & $\leq 5$ & $\leq 30$ & $\leq 100$ & $\leq 255$ \\
    \hline Our method			&  0.39  &  0.93   & 2.63     & 7.21    \\
    \hline NaiveGauss			&  0.59  &  1.65   & 3.56     & 6.9	\\
    \hline AnsGauss			&  0.87  &  2.33   & 4.61     & 8.35	\\
    \hline RL-MRS                       &  0.35  &  1.76   & 4.31     & 9.5	\\
    \hline RL                           &  1.97  &  5.07   & 9.53     & 15.68   \\
    \hline
  \end{tabular}
%   \begin{tabular}{|c|c|c|c|c|}
%     \hline  & \multicolumn{4}{|c|}{Intensity regime}\\
%     \hline \multicolumn{1}{|r|}{Method} & 0 \ldots 5 & 0 \ldots 30 & 0 \ldots 100 & 0 \ldots 255 \\
%     \hline M1                           & 0.25       & 1.87        & 17.67        & 130.42       \\
%     \hline M2                           & 0.61       & 4.72        & 22.59        & 91.41        \\
%     \hline M3                           & 1.26       & 9.03        & 35.71        & 121.08       \\
%     \hline M4                           & 0.24       & 6.79        & 43.67        & 229.76       \\
%     \hline M5                           & 7.5        & 44.29       & 150.12       & 403.65       \\
%     \hline
%   \end{tabular}
  \caption{Mean $\ell_1$-error of all algorithms as a function of the intensity level.}
  \label{tab:intens}
\end{table}

\section{Conclusion}
In this paper, we presented a sparsity-based fast iterative thresholding deconvolution algorithm that takes account of the presence of Poisson noise. A careful theoretical characterization of the algorithm and its solution is provided. The encouraging experimental results clearly confirm the capabilities of our approach.

\bibliographystyle{IEEEbib}
\bibliography{paperbib}

\end{document}